\documentclass[12pt,a4paper]{amsart}
\usepackage{amscd,amsfonts,amsmath,amsthm,amssymb,verbatim,mathrsfs}
\usepackage[dvips]{graphicx}
\usepackage{graphics,hyperref}
\usepackage{enumitem}
\usepackage{psfrag}
\usepackage{latexsym}

\def\frk{\frak}               

\def\Phi{{\frk n}}
\def\Phi{{\frk N}}
\def\opn#1#2{\def#1{\operatorname{#2}}} 
\opn\chara{char} \opn\length{\ell} \opn\pd{pd} \opn\rk{rk}
\opn\projdim{proj\,dim} \opn\injdim{inj\,dim} \opn\rank{rank}
\opn\depth{depth} \opn\grade{grade} \opn\height{height}
\opn\embdim{emb\,dim} \opn\codim{codim} \opn\sgn{sgn}

\opn\Tr{Tr} \opn\bigrank{big\,rank}
\opn\superheight{superheight}\opn\lcm{lcm}
\opn\trdeg{tr\,deg}
\opn\reg{reg} \opn\lreg{lreg} \opn\ini{in} \opn\lpd{lpd}
\opn\size{size}\opn\bigsize{bigsize}
\opn\cosize{cosize}\opn\bigcosize{bigcosize}
\opn\sdepth{sdepth}\opn\sreg{sreg}
\opn\link{link}\opn\fdepth{fdepth}
\opn\div{div} \opn\Div{Div} \opn\cl{cl} \opn\Cl{Cl} \opn\Cor{Cor}
\opn\Spec{Spec} \opn\Supp{Supp} \opn\supp{supp} \opn\Sing{Sing}
\opn\Ass{Ass} \opn\Min{Min}\opn\Mon{Mon} \opn\dstab{dstab} \opn\astab{astab}
\opn\Ann{Ann} \opn\Rad{Rad} \opn\Soc{Soc} \opn\Gr{Gr}
\opn\Im{Im} \opn\Ker{Ker} \opn\Coker{Coker} \opn\Am{Am}
\opn\Hom{Hom} \opn\Tor{Tor} \opn\Ext{Ext} \opn\End{End}
\opn\Aut{Aut} \opn\id{id} \opn\span{span}

\opn\nat{nat}
\opn\pff{pf}
\opn\Pf{Pf} \opn\GL{GL} \opn\SL{SL} \opn\mod{mod} \opn\ord{ord}
\opn\Gin{Gin} \opn\Hilb{Hilb}\opn\sort{sort} \opn\Gale{Gale}
\opn\aff{aff} \opn\conv{conv} \opn\relint{relint} \opn\st{st}   \opn\cone{cone}
\opn\lk{lk} \opn\cn{cn} \opn\core{core} \opn\vol{vol}
\opn\link{link} \opn\star{star}\opn\lex{lex} \opn\Gr{Gr}
\opn\gr{gr}

\def\pot#1#2{#1[\kern-0.28ex[#2]\kern-0.28ex]}

\opn\dirlim{\underrightarrow{\lim}}
\opn\inivlim{\underleftarrow{\lim}}
\def\Implies{\ifmmode\Longrightarrow \else
        \unskip${}\Longrightarrow{}$\ignorespaces\fi}
\def\implies{\ifmmode\Rightarrow \else
        \unskip${}\Rightarrow{}$\ignorespaces\fi}
\def\iff{\ifmmode\Longleftrightarrow \else
        \unskip${}\Longleftrightarrow{}$\ignorespaces\fi}

\let\:=\colon
\newtheorem{Theorem}{Theorem}[section]
\newtheorem{Lemma}[Theorem]{Lemma}
\newtheorem{Corollary}[Theorem]{Corollary}
\newtheorem{Proposition}[Theorem]{Proposition}
\newtheorem{Remark}[Theorem]{Remark}

\newtheorem{Example}[Theorem]{Example}

\let\epsilon\varepsilon
\let\kappa=\varkappa
\makeatletter

\def\qed{\ifhmode\textqed\fi
      \ifmmode\ifinner\quad\qedsymbol\else\dispqed\fi\fi}
\def\textqed{\unskip\nobreak\penalty50
       \hskip2em\hbox{}\nobreak\hfil\qedsymbol
       \parfillskip=0pt \finalhyphendemerits=0}
\def\dispqed{\rlap{\qquad\qedsymbol}}

\opn\dis{dis}
\def\pnt{{\raise0.5mm\hbox{\large\bf.}}}

\opn\Lex{Lex}

\begin{document}


\title{Toric splittings}

\author[1]{Anargyros Katsabekis}
\thanks{Corresponding author: Anargyros Katsabekis}

\author[2]{Apostolos Thoma}

\address{Anargyros Katsabekis, Department of Mathematics, University of Ioannina, Ioannina 45110, Greece}
\email{katsampekis@uoi.gr}

\address{Apostolos Thoma, Department of Mathematics, University of Ioannina, Ioannina 45110, Greece}
\email{athoma@uoi.gr}


\begin{abstract}
 The toric ideal $I_A$ is splittable if it has a toric splitting; namely, if there exist toric ideals $I_{A_1}, I_{A_2}$ such that $I_A=I_{A_1}+I_{A_2}$ and $I_{A_i}\not =I_{A}$ for all $1 \leq i \leq 2$. We provide a necessary and sufficient condition for a toric ideal to be splittable in terms of $A$, and we apply it to prove or disprove that certain classes of toric ideals are splittable.
\end{abstract}
\maketitle

\section{Introduction}

Let $A=\{{\bf a}_{1},\ldots,{\bf a}_{n}\}$ be a vector configuration in $\mathbb{Z}^{m}$ such that the affine semigroup $\mathbb{N}A=\{l_{1}{\bf a}_{1}+\cdots+l_{n}{\bf a}_{n} \mid l_{i} \in \mathbb{N}\}$ is pointed. Let ${\rm ker}_{\mathbb{Q}}(A)=\{{\bf u}=(u_{1},\ldots,u_{n}) \in \mathbb{Q}^{n} \mid u_{1}{\bf a}_{1}+\cdots+u_{n}{\bf a}_{n}={\bf 0}\}$ and ${\rm ker}_{\mathbb{Z}}(A)={\rm ker}_{\mathbb{Q}}(A) \cap \mathbb{Z}^{n}$. Recall that $\mathbb{N}A$ is pointed if ${\rm ker}_{\mathbb{Z}}(A) \cap \mathbb{N}^{n}=\{{\bf 0}\}$. Consider the polynomial ring $K[x_1,\ldots,x_n]$ over a field $K$. The {\em toric ideal} $I_A$ is the kernel of the $K$-algebra homomorphism $\phi: K[x_{1},\ldots,x_{n}] \longrightarrow K[t_{1},\ldots,t_{m},t_{1}^{-1},\ldots,t_{m}^{-1}]$ given by $\phi(x_{i})={\bf t}^{{\bf a}_{i}}$ for all $i=1,\ldots,n$ and it is generated by all the binomials $B({\bf u}):={\bf x}^{{\bf u}_{+}}-{\bf x}^{{\bf u}_{-}}$ such that ${\bf u}_{+}-{\bf u}_{-} \in {\rm ker}_{\mathbb{Z}}(A)$, where ${\bf u}_{+} \in \mathbb{N}^{n}$ and ${\bf u}_{-} \in \mathbb{N}^{n}$ denote the positive and negative part of ${\bf u}={\bf u}_{+}-{\bf u}_{-}$, respectively (see, for instance, \cite[Lemma 4.1]{St}).

A toric ideal $I_A$ is called {\em splittable} if there exist toric ideals $I_{A_1}$ and $I_{A_2}$ such that $I_A=I_{A_1}+I_{A_2}$ and $I_{A_i} \neq I_{A}$ for each $i=1, 2$. In this case, $I_A=I_{A_1}+I_{A_2}$ is called a {\em toric splitting}. The problem of determining when $I_A$ is splittable has attracted the attention of many researchers in recent years; see, for instance, \cite{FHKT, GS, KT2}. The splittings of a toric ideal $I_A$ in general, and especially for the toric ideals $I_G$ of graphs, were considered by G. Favacchio, J. Hofscheier, G. Keiper, and A. Van Tuyl in \cite{FHKT}, in order to identify when $I_A$ admits a toric splitting such that the graded Betti numbers of $I_A$ could be computed using the graded Betti numbers of two other toric ideals, $I_{A_1}$ and $I_{A_2}$.

At the end of Section 5 in \cite{FHKT}, G. Favacchio et al. pose a general question: for which vector configurations $A$ can we identify $A_1$ and $A_2$ such that their corresponding toric ideals satisfy $I_A=I_{A_1}+I_{A_2}$? Can we classify when $I_A$ is a splittable toric ideal in terms of $A$? The problem of determining when the toric ideal $I_G$ of a graph $G$ admits a toric splitting of the form $I_G = I_{G_1} + I_{G_2}$, in terms of the graph $G$, was solved in \cite{KT2}, where $I_{G_i}$ (for $i=1, 2$) denotes the toric ideal of a subgraph $G_i$ of $G$.  

This article seeks to answer the general question posed by G. Favacchio et al. More precisely, in Section 2, we show that $I_A$ is splittable if and only if there exists a minimal system of binomial generators $\{B({\bf u}) \mid {\bf u}\in C\subset {\rm ker}_{\mathbb{Z}}(A)\}$ of the toric ideal $I_A$, and sets $C_{1}$ and $C_2$ such that $C=C_1 \cup C_2$, ${\rm span}_{\mathbb{Q}}(C_1)\subsetneqq {\rm ker}_{\mathbb{Q}}(A)$, and ${\rm span}_{\mathbb{Q}}(C_2)\subsetneqq {\rm ker}_{\mathbb{Q}}(A)$ (see Theorem \ref{main}). We prove that if a toric ideal has a splitting, then it has finitely many different splittings (see Theorem \ref{finitelymany}). 

Section 3 is devoted to the applications of the main result obtained in this paper, namely, Theorem \ref{main}. More precisely, we prove (see Theorems \ref{GluingBasic} and \ref{Gorenstein}) that the ideal $I_A$ is splittable in the following cases, and many other related to them,: \begin{enumerate} \item The semigroup $\mathbb{N}A$ is obtained by gluing the semigroups $\mathbb{N}A_1$ and $\mathbb{N}A_2$.
\item $A=\{a_{1},\ldots,a_{4}\}$ is a set of positive integers such that the corresponding monomial curve, with parametrization $x_{i}=t^{a_i}$ for $1 \leq i \leq 4$, is Gorenstein.
\end{enumerate}

Finally, we show that the toric ideal of a complete bipartite graph is not splittable (see Theorem \ref{CompleteBipartite}).\\

\section{Splitting Criterion for $I_A$}
\label{section}

Let $A=\{{\bf a}_{1},\ldots,{\bf a}_{n}\}$ be a vector configuration in $\mathbb{Z}^{m}$. Given a set $C=\{{\bf c}_{1},\ldots,{\bf c}_{k}\}\subset \mathbb{Z}^n$, we shall denote by ${\rm span}_{\mathbb{Q}}(C)$ the $\mathbb{Q}$-vector space generated by the vectors of $C$, namely, $${\rm span}_{\mathbb{Q}}(C)=\{\lambda_{1}{\bf c}_{1}+\cdots+\lambda_{k}{\bf c}_{k} \mid \lambda_{i} \in \mathbb{Q}\}.$$

\begin{Theorem} \label{main} The toric ideal $I_A$ is splittable if and only if there exists a minimal system of binomial generators $\{B({\bf u}) \mid {\bf u}\in C\subset {\rm ker}_{\mathbb{Z}}(A)\}$ of the toric ideal $I_A$, and sets $C_{1}$ and $C_2$ such that $C=C_1 \cup C_2$, ${\rm span}_{\mathbb{Q}}(C_1)\subsetneqq {\rm ker}_{\mathbb{Q}}(A)$, and ${\rm span}_{\mathbb{Q}}(C_2)\subsetneqq {\rm ker}_{\mathbb{Q}}(A)$.   
\end{Theorem}
{\em \noindent Proof.} Suppose that $I_{A_1}+I_{A_2}$ is a splitting of $I_A$. Thus, ${\rm ker}_{\mathbb{Z}}(A_1)\subsetneqq {\rm ker}_{\mathbb{Z}}(A)$ and ${\rm ker}_{\mathbb{Z}}(A_2)\subsetneqq {\rm ker}_{\mathbb{Z}}(A)$.
Let $\{B({\bf v}_{1}),\ldots,B({\bf v}_{r})\}$ and $\{B({\bf w}_{1}),\ldots,B({\bf w}_{t})\}$ be minimal systems of binomial generators of $I_{A_1}$ and $I_{A_2}$, respectively. Then, since $I_A=I_{A_1}+I_{A_2}$ and $I_A$ is positively graded, there exists a minimal system of binomial generators of $I_A$ in the form $\{B({\bf u}) \mid {\bf u} \in C\subset \{{\bf v}_{1},\ldots, {\bf v}_{r}, {\bf w}_{1},\ldots, {\bf w}_{t}\}\}$. This means that each minimal generator of $I_A$ is either a minimal generator of $I_{A_1}$ or a minimal generator of $I_{A_2}$. Let $C_1=C\cap \{{\bf v}_{1},\ldots,{\bf v}_{r}\}$ and $C_2=C\cap \{{\bf w}_{1},\ldots,{\bf w}_{t}\}$. Thus $C=C_1 \cup C_2$. 
We claim that ${\rm span}_{\mathbb{Q}}(C_1)\subsetneqq {\rm ker}_{\mathbb{Q}}(A)$ and ${\rm span}_{\mathbb{Q}}(C_2)\subsetneqq {\rm ker}_{\mathbb{Q}}(A)$. Suppose not. Then ${\rm span}_{\mathbb{Q}}(C_i)= {\rm ker}_{\mathbb{Q}}(A)$ for some $i$ such that $1 \leq i \leq 2$. Thus, ${\rm span}_{\mathbb{Q}}(C_i)\cap \mathbb{Z}^n = {\rm ker}_{\mathbb{Q}}(A)\cap \mathbb{Z}^n={\rm ker}_{\mathbb{Z}}(A)$. 
Note also that $C_i\subset {\rm ker}_{\mathbb{Z}}(A_i)\subsetneqq {\rm ker}_{\mathbb{Z}}(A)$. Thus, ${\rm span}_{\mathbb{Q}}(C_i)\subset {\rm ker}_{\mathbb{Q}}(A_i)\subsetneqq {\rm ker}_{\mathbb{Q}}(A)$. Consequently, $I_A \subset I_{A_i}\subsetneqq I_A$, which is a contradiction. Thus, ${\rm span}_{\mathbb{Q}}(C_1)\subsetneqq {\rm ker}_{\mathbb{Q}}(A)$ and ${\rm span}_{\mathbb{Q}}(C_2)\subsetneqq {\rm ker}_{\mathbb{Q}}(A)$.

Conversely, suppose that there exists a minimal system of binomial generators $\{B({\bf u}) \mid {\bf u}\in C\subset {\rm ker}_{\mathbb{Z}}(A)\}$ and sets $C_{1}$ and $C_2$ such that $C=C_1 \cup C_2$, and ${\rm span}_{\mathbb{Q}}(C_1)\subsetneqq {\rm ker}_{\mathbb{Q}}(A)$ and ${\rm span}_{\mathbb{Q}}(C_2)\subsetneqq {\rm ker}_{\mathbb{Q}}(A)$. The set ${\rm span}_{\mathbb{Q}}(C_1)$ is a subspace of $\mathbb{Q}^{n}$. Consider any set of integer vectors $A_1$  that defines the orthogonal complement of the vector space ${\rm span}_{\mathbb{Q}}(C_1)$. Consequently, ${\rm ker}_{\mathbb{Q}}(A_{1})={\rm span}_{\mathbb{Q}}(C_1)$. Similarly, consider a set of integer vectors $A_2$ such that ${\rm ker}_{\mathbb{Q}}(A_{2})={\rm span}_{\mathbb{Q}}(C_2)$.
 
We claim that $I_A=I_{A_1}+I_{A_2}$  is a splitting of $I_A$.
Since $C_i\subset {\rm ker}_{\mathbb{Z}}(A)$, it follows that
${\rm ker}_{\mathbb{Z}}(A_{i}) \subset {\rm ker}_{\mathbb{Z}}(A)$. Therefore, $I_{A_1}\subset I_A$ and $I_{A_2}\subset I_A$, and thus $I_{A_1}+I_{A_2} \subset I_A$. For each minimal generator $B({\bf u}_{i})$, where $1 \leq i \leq s$, of $I_A$, we have that ${\bf u}_{i}\in C_1$ or ${\bf u}_{i} \in C_2$. If ${\bf u}_{i} \in C_{1}$, then ${\bf u}_{i} \in {\rm ker}_{\mathbb{Q}}(A_{1})$ and thus $B({\bf u}_{i})\in I_{A_1}$. Similarly, if ${\bf u}_{i} \in C_2$, we conclude that $B({\bf u}_{i})\in I_{A_2}$. Therefore, every minimal generator of $I_A$ is in $I_{A_1}$ or $I_{A_2}$. We conclude that $I_A\subset I_{A_1}+I_{A_2}$, and thus $I_A=I_{A_1}+I_{A_2}$. 

Note also that $I_{A_i}\not =I_A$ for $1 \leq i \leq 2$. Otherwise, if $I_{A_i} =I_A$, then ${\rm ker}_{\mathbb{Z}}(A)={\rm ker}_{\mathbb{Z}}(A_{i})$, and therefore ${\rm ker}_{\mathbb{Q}}(A)={\rm span}_{\mathbb{Q}}(C_i)$, which is a contradiction.
\hfill $\square$

\begin{Example} \label{many} {\rm  Let $A=\{4, 6, 11, 13\}$ be a set of integers. The toric ideal $I_A$ is minimally generated by the binomials $B({\bf u}_{i})$ for $1 \leq i \leq 6$, where ${\bf u}_1=(3,-2,0,0)$, ${\bf u}_2=(-4,-1,2,0)$, ${\bf u}_3=(-1,0,-2,2)$, ${\bf u}_4=(1,-1,-1,1)$, ${\bf u}_5=(2,-1,1,-1)$ and ${\bf u}_6=(6,0,-1,-1)$. Therefore, from \cite[Proposition 8.3.1]{Vil}, ${\rm ker}_{\mathbb{Z}}(A)=\langle {\bf u}_{i} \mid  1 \leq i \leq 6 \rangle$, so ${\rm ker}_{\mathbb{Q}}(A)={\rm span}_{\mathbb{Q}}({\bf u}_{i} \mid  1 \leq i \leq 6)$. Since ${\bf u}_4=\frac{{\bf u}_1}{2}+\frac{{\bf u}_3}{2}, {\bf u}_5=\frac{{\bf u}_1}{2}-\frac{{\bf u}_3}{2}$ and ${\bf u}_6=\frac{{\bf u}_1}{2}-{\bf u}_2-\frac{{\bf u}_3}{2}$, it follows that ${\rm ker}_{\mathbb{Q}}(A)={\rm span}_{\mathbb{Q}}({\bf u}_{1},{\bf u}_{2}, {\bf u}_{3})$. Furthermore, since $\{{\bf u}_{1},{\bf u}_{2}, {\bf u}_{3}\}$ is linearly independent, the set $\{{\bf u}_1, {\bf u}_2, {\bf u}_3\}$ is a basis of the vector space ${\rm ker}_{\mathbb{Q}}(A)$. Let $C=\{{\bf u}_i \mid 1 \leq i \leq 6\}$, and suppose that there are sets $C_{1}$ and $C_2$ such that $C=C_1 \cup C_2$, ${\rm span}_{\mathbb{Q}}(C_1)\subsetneqq {\rm ker}_{\mathbb{Q}}(A)$, and ${\rm span}_{\mathbb{Q}}(C_2)\subsetneqq {\rm ker}_{\mathbb{Q}}(A)$.  Since the set $\{{\bf u}_1, {\bf u}_2, {\bf u}_3\}$ is a basis of the vector space ${\rm ker}_{\mathbb{Q}}(A)$, exactly two of the vectors ${\bf u}_1, {\bf u}_2, {\bf u}_3 $ belong to one of the $C_i$, and the remaining one belongs to $C_j$, where $\{i, j\}=\{1, 2\}$. We may assume that two of them belong to $C_1$, and the remaining one belongs to $C_2$. There are three cases.
\begin{enumerate}
\item ${\bf u}_1, {\bf u}_2\in C_1$ and ${\bf u}_3 \in C_2$. The sets $\{{\bf u}_1, {\bf u}_2, {\bf u}_4\}$, $\{{\bf u}_1, {\bf u}_2, {\bf u}_5\}$, and $\{{\bf u}_1, {\bf u}_2, {\bf u}_6\}$ are bases of ${\rm ker}_{\mathbb{Q}}(A)$. Therefore, ${\bf u}_4, {\bf u}_5, {\bf u}_6$ belong to the set $C_2$. However, the set
$C_2$ generates ${\rm ker}_{\mathbb{Q}}(A)$; thus, ${\rm span}_{\mathbb{Q}}(C_2)= {\rm ker}_{\mathbb{Q}}(A)$, which is a contradiction.
\item ${\bf u}_2, {\bf u}_3\in C_1$, and ${\bf u}_1 \in C_2$. The sets $\{{\bf u}_2, {\bf u}_3, {\bf u}_4\}$, $\{{\bf u}_2, {\bf u}_3, {\bf u}_5\}$, and $\{{\bf u}_2, {\bf u}_3, {\bf u}_6\}$ are bases of ${\rm ker}_{\mathbb{Q}}(A)$. Therefore, 
${\bf u}_4, {\bf u}_5, {\bf u}_6$ belong to the set $C_2$. But then the set $C_2$ generates ${\rm ker}_{\mathbb{Q}}(A)$, and thus ${\rm span}_{\mathbb{Q}}(C_2)= {\rm ker}_{\mathbb{Q}}(A)$, which is a contradiction.
\item ${\bf u}_1$ and ${\bf u}_3$ are in $C_1$, and ${\bf u}_2$ is in $C_2$. The set $\{{\bf u}_1, {\bf u}_3, {\bf u}_6\}$
is a basis of ${\rm ker}_{\mathbb{Q}}(A)$, so ${\bf u}_6$ belongs to the set $C_2$. The set $\{{\bf u}_2, {\bf u}_6, {\bf u}_4\}$
is a basis of ${\rm ker}_{\mathbb{Q}}(A)$; therefore, ${\bf u}_4$
belongs to the set $C_1$. Thus, in this case, we may have $C_1=\{{\bf u}_1, {\bf u}_3, {\bf u}_4, {\bf u}_5\}$ and $C_2=\{{\bf u}_2, {\bf u}_5, {\bf u}_6\}$; or $C_1=\{{\bf u}_1, {\bf u}_3, {\bf u}_4\}$ and $C_2=\{{\bf u}_2, {\bf u}_5, {\bf u}_6\}$; or $C_1=\{{\bf u}_1, {\bf u}_3, {\bf u}_4, {\bf u}_5\}$ and $C_2=\{{\bf u}_2, {\bf u}_6 \}$. In all three cases, ${\rm span}_{\mathbb{Q}}(C_{1})={\rm span}_{\mathbb{Q}}({\bf u}_1, {\bf u}_3)$, since ${\bf u}_{4}=\frac{1}{2}{\bf u}_{1}+\frac{1}{2}{\bf u}_{3}$ and ${\bf u}_{5}=\frac{1}{2}{\bf u}_{1}-\frac{1}{2}{\bf u}_{3}$. Similarly, ${\rm span}_{\mathbb{Q}}(C_{2})={\rm span}_{\mathbb{Q}}({\bf u}_2, {\bf u}_6)$, because ${\bf u}_5={\bf u}_2+{\bf u}_6$. Also, ${\rm span}_{\mathbb{Q}}({\bf u}_1, {\bf u}_3) \neq {\rm ker}_{\mathbb{Q}}(A)$ and ${\rm span}_{\mathbb{Q}}({\bf u}_2, {\bf u}_6) \neq {\rm ker}_{\mathbb{Q}}(A)$, since ${\rm ker}_{\mathbb{Q}}(A)$ is three dimensional. By Theorem \ref{main}, the ideal $I_A$ is splittable. Let $A_{1}$ be the set of columns of the matrix
 \[
\left( \begin{array}{cccc}
 2 &  3 &  3 & 4  \\
 4 & 6 & 11 & 13
\end{array} \right)
	 \in \mathbb{Z}^{2 \times 4}
\] and $A_{2}$  be the set of columns of the matrix \[
\left( \begin{array}{cccc}
1 &  8 &  6 & 0  \\
0 & -2 & -1 & 1
\end{array} \right)
	 \in \mathbb{Z}^{2 \times 4}.
\]
Then ${\rm ker}_{\mathbb{Q}}(A_1)={\rm span}_{\mathbb{Q}}(C_1)$ and ${\rm ker}_{\mathbb{Q}}(A_2)={\rm span}_{\mathbb{Q}}(C_2)$. From the proof of Theorem \ref{main}, we have that 
$I_{A}=I_{A_1}+I_{A_2}$ is a splitting of $I_A$.
\end{enumerate}

Note that, from the analysis above, it seems there exists a unique splitting of $I_A$. This is not true; there are several different splittings of $I_A$. The reason is that the analysis above was based on a given minimal system of binomial generators of $I_A$, and in this example, we have 18 different minimal systems of binomial generators. Note that, from the binomials $B({\bf u}_i), 1 \leq i \leq 6$, only $B({\bf u}_1), B({\bf u}_4),$ and $B({\bf u}_5)$ are indispensable, see \cite{ChKTh}. Recall that a binomial $B \in I_A$ is indispensable if every system of binomial generators of $I_A$ contains $B$ or $-B$. The binomial $B({\bf u}_2)$ can be replaced in a minimal system of generators by $B({\bf u}'_2)$ for ${\bf u}'_2=(-1, -3,2,0)={\bf u}_1+{\bf u}_2$, $B({\bf u}_3)$ by $B({\bf u}'_3)$ or $B({\bf u}''_3)$ for ${\bf u}'_3=(-5, -1, 0, 2)$ and ${\bf u}''_3=(-2, -3, 0, 2)$, and $B({\bf u}_6)$ by $B({\bf u}'_6)$ or $B({\bf u}''_6)$ for ${\bf u}'_6=(3,2, -1, -1)$ and ${\bf u}''_6=(0, 4, -1, -1)$. For example, we can choose the minimal system $\{B({\bf u}) \mid {\bf u}\in C'\subset {\rm ker}_{\mathbb{Z}}(A)\}$ of binomial generators of $I_A$,  where $C'=\{{\bf u}_1, {\bf u}'_2, {\bf u}'_3, {\bf u}_4, {\bf u}_5, {\bf u}'_6\}$. Considering again the different possible cases, we can see that the only way to have a splitting of $I_A$ is by $C'_{1}=\{{\bf u}_{1}, {\bf u}_{4}\}$ and $C'_{2}=\{{\bf u}'_{2}, {\bf u}'_3, {\bf u}_{5}, {\bf u}'_{6}\}$. We have that ${\rm span}_{\mathbb{Q}}(C'_{1})\not={\rm ker}_{\mathbb{Q}}(A)$ and ${\rm span}_{\mathbb{Q}}(C'_{2})\not={\rm ker}_{\mathbb{Q}}(A)$. Let $A'_{1}$ be the set of columns of the matrix
 \[
\left( \begin{array}{cccc}
 4 &  6 &  11 & 13  \\
 0 & 0 & 1 & 1
\end{array} \right)
	 \in \mathbb{Z}^{2 \times 4}
\] and $A'_{2}$  be the set of columns of the matrix \[
\left( \begin{array}{cccc}
1 &  -5 &  -7 & 0  \\
0 & 2 & 3 & 1
\end{array} \right)
	 \in \mathbb{Z}^{2 \times 4}.
\]
Then ${\rm ker}_{\mathbb{Q}}(A'_1)={\rm span}_{\mathbb{Q}}(C'_1)$ and ${\rm ker}_{\mathbb{Q}}(A'_2)={\rm span}_{\mathbb{Q}}(C'_2)$. From the proof of Theorem \ref{main}, we have that 
$I_{A}=I_{A'_1}+I_{A'_2}$ is a splitting of $I_A$. This is a different splitting from $I_{A}=I_{A_1}+I_{A_2}$, since ${\rm ker}_{\mathbb{Q}}(A_1) = {\rm ker}_{\mathbb{Q}}(A'_1)$ but ${\rm ker}_{\mathbb{Q}}(A_2)\not ={\rm ker}_{\mathbb{Q}}(A'_2)$.
Note that there are infinite sets of vectors that define the same toric ideal; all of them have the same kernel.
}\end{Example}

\begin{Example} \label{none} {\rm  Let $A=\{20, 24, 25, 31\}$ be a set of integers. The toric ideal $I_A$ is generic and thus has a unique minimal system of binomial generators (see \cite[Remark 4.4.3]{PS}). There are seven minimal generators of $I_A$ in the form $B({\bf u})$, where ${\bf u}$ is one of the following vectors: $(4, -1, -1, -1)$, $(-2, 4, -1, -1)$, $(-1, -1, 3, -1)$, $(-1, -2, -1, 3)$, $(2,3,-2,-2)$, $(3, -2, 2, -2)$, and $(3, -3, -2, 2)$. We denote the first four vectors as ${\bf a}, {\bf b}, {\bf c}, {\bf d}$, in no particular order. Then we have ${\bf a}+{\bf b}+{\bf c}+{\bf d}={\bf 0}$, and any three of them generate the vector space ${\rm ker}_{\mathbb{Q}}(A)$. The remaining three vectors, up to a sign change, are ${\bf a}+{\bf b}=-{\bf c}-{\bf d}$, ${\bf a}+{\bf c}=-{\bf b}-{\bf d}$, and ${\bf a}+{\bf d}=-{\bf b}-{\bf c}$. Suppose that $I_A=I_{A_1}+I_{A_2}$ is a splitting. Then, by Theorem \ref{main}, there exist sets $C_{1}$ and $C_2$ such that $C=\{{\bf a}, {\bf b}, {\bf c}, {\bf d}, {\bf a}+{\bf b}, {\bf a}+{\bf c}, {\bf a}+{\bf d}\}=C_1 \cup C_2$, with ${\rm span}_{\mathbb{Q}}(C_1)\subsetneqq {\rm ker}_{\mathbb{Q}}(A)$ and ${\rm span}_{\mathbb{Q}}(C_2)\subsetneqq {\rm ker}_{\mathbb{Q}}(A)$. Since any three of ${\bf a}$, ${\bf b}$, ${\bf c}$, and ${\bf d}$ generate ${\rm ker}_{\mathbb{Q}}(A)$, the set $C_1$
should contain exactly two of them, say ${\bf a}$ and ${\bf b}$, and the set $C_2$ should contain the other two, ${\bf c}$ and ${\bf d}$. Since ${\bf a}, {\bf b}$, and ${\bf a}+{\bf c}$ generate ${\rm ker}_{\mathbb{Q}}(A)$, the set $C_2$ should also contain the vector ${\bf a}+{\bf c}$. However, the vectors ${\bf c}, {\bf d}$, and $ {\bf a}+{\bf c}$ in $C_2$ generate ${\rm ker}_{\mathbb{Q}}(A)$, which contradicts Theorem \ref{main}. Thus, the toric ideal $I_A$ does not have a splitting.}
    
\end{Example}
\begin{Remark} {\em In Examples \ref{many} and \ref{none}, we observed that a toric ideal may have many splittings or none at all. The number of different splittings can be huge. In \cite{KT2}, we studied subgraph splittings of the toric ideal $I_G$ of a graph $G$, see Section 3 for the definition of this ideal. More precisely, \cite[Proposition 3.9]{KT2} shows that if $n \geq 4$ and $I_{K_n}=I_{G_1}+I_{G_2}$ is a splitting, where $G_1$ and $G_2$ are subgraphs of the complete graph $K_n$ on $n$ vertices, then there exists a cycle $w=(a,b,c,d)$ in $K_n$ such that 
\begin{enumerate}
    \item $G_1={K_n}\symbol{92} a$ and $G_2={K_n}\symbol{92} b$ or
    \item $G_1={K_n}\symbol{92} \{a,c\}$ and $G_2={K_n}\symbol{92} b$ or 
    \item $G_1={K_n}\symbol{92} a$ and $G_2={K_n}\symbol{92} \{b,d\}$ or
    \item $G_1={K_n}\symbol{92} \{a,c\}$ and $G_2={K_n}\symbol{92} \{b,d\}$.
\end{enumerate} 
There are $3\binom{n}{3}$ splittings of type (1), $12\binom{n}{4}$ of types (2) and (3), and $3\binom{n}{4}$ of type (4). Thus, the toric ideal of $K_n$ has exactly $3\binom{n}{3}+15\binom{n}{4}$ splittings into toric ideals of subgraphs, so $I_{K_n}$ has at least $3\binom{n}{3}+15\binom{n}{4}$ splittings. For large values of $n$, this represents a huge number.}
\end{Remark}

The following theorem states that if a toric ideal $I_A$ has a splitting, then it has finitely many different toric splittings:

\begin{Theorem} \label{finitelymany} Let $I_A$ be a splittable toric ideal, then it has finitely many splittings.
\end{Theorem}
{\em \noindent Proof.}
Since ${\rm ker}_{\mathbb{Z}}(A) \cap \mathbb{N}^{n}=\{{\bf 0}\}$, the toric ideal has finitely many distinct minimal systems of binomial generators. See \cite[Theorem 2.9]{ChKTh} for a formula providing the number of distinct minimal systems of binomial generators of a toric ideal $I_A$. For each minimal system of binomial generators $\{B({\bf u}) \mid {\bf u}\in C\subset {\rm ker}_{\mathbb{Z}}(A)\}$ of a toric ideal $I_A$, there are finitely many ways to express $C$ as $C_1 \cup C_2$. Therefore, there are finitely many toric ideals $I_{A_1}$ and $I_{A_2}$ such that $I_A=I_{A_1}+I_{A_2}$.
\hfill $\square$

\section{Applications}

In this section, we first study the case of the gluing of two semigroups.
Let $A \subset \mathbb{Z}^{n}$ be a finite set of vectors. Let $A_1=\{{\bf a}_{11},\ldots, {\bf a}_{1s}\}$ and $A_2=\{{\bf a}_{21},\ldots , {\bf a}_{2t}\}$ be two non-empty subsets of $A$ such that $A=A_1\cup A_2$ and $A_1\cap A_2=\emptyset$. We say that the semigroup $\mathbb{N}A$ is obtained by {\em gluing} the semigroups $\mathbb{N}A_1$ and $\mathbb{N}A_2$ if there exists a nonzero ${\bf a}\in \mathbb{N}A_1\cap \mathbb{N}A_2$, such that $\mathbb{Z}A_1\cap \mathbb{Z}A_2=\mathbb{Z}{\bf a}$. If $\mathbb{N}A$ is obtained by gluing the semigroups $\mathbb{N}A_1$ and $\mathbb{N}A_2$, then ${\bf a}$ can be written as ${\bf a}=c_1{\bf a}_{11}+\cdots +c_s{\bf a}_{1s}=d_1{\bf a}_{21}+\cdots +d_t{\bf a}_{2t}$, where all $c_i\in \mathbb{N}, d_i\in \mathbb{N}$. In terms of ideals, the gluing is equivalent to $I_A=I_{A_1}+I_{A_2}+<x_{11}^{c_1}\cdots x_{1s}^{c_s}-x_{21}^{d_1}\cdots x_{2t}^{d_t}>$; see \cite{Ros}.

  In \cite[Theorem 2.2]{GS},  P. Gimenez and H. Srinivasan  showed  that if the semigroup $\mathbb{N}A$ is obtained by a particular gluing of two homogeneous semigroups, then $I_A$ is splittable. 
The next theorem generalizes Theorem 2.2 of \cite{GS}.

\begin{Theorem} \label{GluingBasic} 
Let $\mathbb{N}A$ be the semigroup obtained by gluing the non-zero semigroups $\mathbb{N}A_1$ and $\mathbb{N}A_2$. Then, the toric ideal $I_A$ is splittable.

\end{Theorem}
{\em \noindent Proof.}
The semigroup  $\mathbb{N}A$ is the gluing of the semigroups $\mathbb{N}A_1$, $\mathbb{N}A_2$;
$I_A=I_{A_1}+I_{A_2}+<x_{11}^{c_1}\cdots x_{1s}^{c_s}-x_{21}^{d_1}\cdots x_{2t}^{d_t}>.$
Let $\{B({\bf u}_{1}),\ldots,B({\bf u}_{k})\}$ and $\{B({\bf v}_{1}),\ldots,B({\bf v}_{l})\}$ be minimal generating sets of $I_{A_1}$ and $I_{A_2}$, respectively. Let ${\bf u}'_{i}=({\bf u}_{i}, {\bf 0}) \in \mathbb{Z}^{s+t}$, $1 \leq i \leq k$, ${\bf v}'_{j}=({\bf 0}, {\bf v}_{j}) \in \mathbb{Z}^{s+t}$, $1 \leq j \leq l$, and ${\bf w}=(c_{1},\ldots,c_{s},-d_{1},\ldots,-d_{t})$, where at least one of $d_{i}$ is nonzero. Then $C=\{{\bf u}'_{1},\ldots,{\bf u}'_{k}, {\bf v}'_{1},\ldots,{\bf v}'_{l}, {\bf w}\}$. Let $C_{1}=\{{\bf u}'_{1},\ldots,{\bf u}'_{k}\}$ and $C_{2}=\{{\bf v}'_{1},\ldots,{\bf v}'_{l}, {\bf w}\}$. Assume that ${\rm span}_{\mathbb{Q}}(C_{1})= {\rm ker}_{\mathbb{Q}}(A)$. Then ${\bf w} \in {\rm ker}_{\mathbb{Q}}(A)$, so ${\bf w}=\lambda_{1} {\bf u}'_{1}+\cdots+\lambda_{k} {\bf u}'_{k}$ for $\lambda_{i} \in \mathbb{Q}$. Thus, $d_{i}=0$ for every $1 \leq i \leq t$, a contradiction. Therefore, ${\rm span}_{\mathbb{Q}}(C_{1}) \subsetneqq {\rm ker}_{\mathbb{Q}}(A)$. Assume that ${\rm span}_{\mathbb{Q}}(C_{2})= {\rm ker}_{\mathbb{Q}}(A)$. Then ${\bf u}'_{i} \in {\rm ker}_{\mathbb{Q}}(A)$, so ${\bf u}'_{i}=\mu_{1} {\bf v}'_{1}+\cdots+\mu_{l} {\bf v}'_{l}+\kappa {\bf w}$ for $\kappa, \mu_{j} \in \mathbb{Q}$. Thus, ${\bf u}_{i}=\kappa {\bf z}$, where
${\bf u}_{i}\in \mathbb{Z}^{s}$ and ${\bf z}=(c_{1},\ldots,c_{s}) \in \mathbb{N}^{s}$, so either ${\bf u}_{i} \in \mathbb{N}^{s}$ or $-{\bf u}_{i} \in \mathbb{N}^{s}$. Note that ${\bf u}_{i}\not ={\bf 0}$ since $B({\bf u}_{i})$ is a minimal generator. Thus, either ${\bf u}'_{i} \in \mathbb{N}^{s+t}$ or $-{\bf u}'_{i} \in \mathbb{N}^{s+t}$, and ${\bf u}'_{i}\not ={\bf 0}$, which contradicts the fact that $\mathbb{N}A$ is pointed. By Theorem \ref{main}, the ideal $I_A$ is splittable. 
\hfill $\square$

\begin{Example} {\rm Let $A=\{85,102,77,88,99\}$, and let $B=\{85,102\}$ and $C=\{77,88,99\}$. Since $A=B \cup C$,  $\mathbb{Z}B \cap \mathbb{Z}C=\mathbb{Z} \{187\}$ and $187\in \mathbb{N}B \cap \mathbb{N}C$,  $\mathbb{N}A$ is obtained by gluing $\mathbb{N}B$ and $\mathbb{N}C$ together. Then $I_A$ is minimally generated by the binomials $x_{1}^6-x_{2}^5$, $x_{3}^5-x_{4}x_{5}^3$, $x_{4}^2-x_{3}x_{5}$, $x_{5}^4-x_{3}^4x_{4}$, $x_{1}x_{2}-x_{4}x_{5}$. Let $A_{1}$ and $A_2$ be the sets consisting of the columns of the $3 \times 5$ matrices \[
\left( \begin{array}{ccccc}
 5 &  6 &  0 & 0 & 11  \\
 0 & 0 & 1 & 0 & 0 \\
 0 & 0 & 0 & 1 & -1
\end{array} \right) \ \textrm{and} \
\left( \begin{array}{ccccc}
 1 &  0 &  0 & 0 & 0  \\
 0 & 1 & 0 & 0 & 0 \\
 0 & 0 & 7 & 8 & 9
\end{array} \right), \ \textrm{respectively}.
\]
Then $I_{A_1}$ is minimally generated by the binomials $x_{1}^6-x_{2}^5$ and $x_{1}x_{2}-x_{4}x_{5}$, and $I_{A_2}$ is minimally generated by the binomials $x_{3}^5-x_{4}x_{5}^3$, $x_{4}^2-x_{3}x_{5}$, and $x_{5}^4-x_{3}^4x_{4}$. So $I_{A}=I_{A_1}+I_{A_2}$ is a splitting of $I_A$.}
\end{Example}

\begin{Corollary}\label{CompleteIntersection} 
Let  $I_A$ be a toric ideal of height greater than or equal to two. If $I_A$ is a complete intersection, then it is splittable.
\end{Corollary}

{\em \noindent Proof.} Since $I_A$ is a complete intersection ideal, we have from \cite[Theorem 3.1]{FMS} that there are two subsets $A_1$ and $A_2$ of $A$ such that $\mathbb{N}A$ is the gluing of $\mathbb{N}A_{1}$ and $\mathbb{N}A_{2}$. By Theorem \ref{GluingBasic}, the ideal $I_A$ is splittable.
\hfill $\square$

\begin{Theorem} \label{height2}
Let $I_A$ be a toric ideal of height 2. Then $I_A$ is splittable if and only if $I_A$ is a complete intersection. 
\end{Theorem}

{\em \noindent Proof.} Suppose first that $I_A$ is a complete intersection. Then, $I_A$ is splittable by Corollary \ref{CompleteIntersection}.\\ Conversely, suppose that $I_A$ is splittable. From Theorem \ref{main}, there exists a minimal system of binomial generators of the toric ideal $I_A$, 
 $S=\{B({\bf w}) \mid {\bf w} \in C\}$, and sets $C_{1}$ and $C_2$ such that $C=C_1 \cup C_2$ with ${\rm span}_{\mathbb{Q}}(C_{1}) \subsetneqq {\rm ker}_{\mathbb{Q}}(A)$ and ${\rm span}_{\mathbb{Q}}(C_{2}) \subsetneqq {\rm ker}_{\mathbb{Q}}(A)$. Then, both the vector spaces ${\rm span}_{\mathbb{Q}}(C_{1})$ and ${\rm span}_{\mathbb{Q}}(C_{2})$ are one-dimensional. We claim that each $C_i$, where $1 \leq i \leq 2$, is a singleton. Suppose that $C_i$ contains at least two elements, ${\bf u}$ and ${\bf w}$. Then there exist  integers $\lambda \neq \mu$ such that $\lambda {\bf w}=\mu {\bf u}$. Since $B({\bf w})$,$B({\bf u})$ are minimal generators of $I_A$, they are irreducible. Thus, the greatest common divisor of their respective coordinates is 1, which implies that $\lambda = \mu$, a contradiction. Therefore, each $C_i$ is a singleton; let $C_{1}=\{{\bf u}\}$ and $C_{2}=\{{\bf v}\}$. Then $S=\{B({\bf u}), B({\bf v})\}$ is a minimal generating set of $I_A$, implying that $I_A$ is a complete intersection.
\hfill $\square$\\

The next theorem shows that the corresponding statement of Theorem \ref{height2} is not true when the height of $I_A$ is greater than two: there are toric ideals that are splittable, but they are not complete intersections.  
\begin{Theorem} \label{BasicMinimal}
If $I_A$ is a toric ideal of height ${\rm ht}(I_A)=s>1$ such that its minimal number of generators is less than or equal to $2s-2$, then $I_A$ is splittable. 
\end{Theorem}

{\em \noindent Proof.} Let $\{B({\bf u}_{1}),\ldots,B({\bf u}_{n})\}$ be a minimal generating set of $I_A$, where $n \leq 2s-2$, and let $C=\{{\bf u}_{1},\ldots,{\bf u}_{n}\}$. Let $C_{1}=\{{\bf u}_{1},\ldots,{\bf u}_{k}\}$ and $C_{2}=\{{\bf u}_{k+1},\ldots,{\bf u}_{n}\}$, where $k=\lfloor \frac{n}{2} \rfloor$. Then, $C=C_{1} \cup C_{2}$. Here, $\lfloor y \rfloor$ denotes the greatest integer less than or equal to $y$. Since the dimension of ${\rm span}_{\mathbb{Q}}(C_{1})$ is less than or equal to $k$ and $k \leq \frac{n}{2}$, we conclude that the dimension of ${\rm span}_{\mathbb{Q}}(C_{1})$ is less than or equal to $\frac{n}{2}$. Therefore, the dimension of ${\rm span}_{\mathbb{Q}}(C_{1})$ is less than or equal to $s-1$ because $n \leq 2s-2$. But the dimension of ${\rm ker}_{\mathbb{Q}}(A)$ is equal to $s$, so the dimension of ${\rm span}_{\mathbb{Q}}(C_{1})$ is strictly smaller than the dimension of ${\rm ker}_{\mathbb{Q}}(A)$. Thus, ${\rm span}_{\mathbb{Q}}(C_{1}) \subsetneqq {\rm ker}_{\mathbb{Q}}(A)$. Note that $\frac{n}{2}<\lfloor \frac{n}{2} \rfloor+1$, and therefore, $\frac{n}{2}+1>n-\lfloor \frac{n}{2} \rfloor$. But $2s-2 \geq n$, so $s \geq \frac{n}{2}+1$ and therefore $s>n-\lfloor \frac{n}{2} \rfloor$. Since the dimension of ${\rm span}_{\mathbb{Q}}(C_{2})$ is less than or equal to $n-\lfloor \frac{n}{2} \rfloor$ and the dimension of ${\rm ker}_{\mathbb{Q}}(A)$ is equal to $s$, we have that the dimension of ${\rm span}_{\mathbb{Q}}(C_{2})$ is strictly smaller than the dimension of ${\rm ker}_{\mathbb{Q}}(A)$. Consequently, ${\rm span}_{\mathbb{Q}}(C_{2}) \subsetneqq {\rm ker}_{\mathbb{Q}}(A)$. By Theorem \ref{main}, the ideal $I_A$ is splittable.
\hfill $\square$\\

The ideal $I_A$ is called an {\em almost complete intersection} if it is minimally generated by ${\rm ht}(I_A)+1$ binomials.

\begin{Corollary}
If $I_A$ is an almost complete intersection ideal of height $s \geq 3$, then $I_A$ is splittable. 
\end{Corollary}

{\em \noindent Proof.} The minimal number of generators of $I_A$ is $s+1$, so $s+1 \leq 2s-2$ since $s \geq 3$. By Theorem \ref{BasicMinimal}, the ideal $I_A$ is splittable.
\hfill $\square$
\bigskip
\newline

The next theorem provides an example showing that the hypotheses of Theorem \ref{BasicMinimal} are not true (the number of minimal generators is greater than $2s-2$), but the corresponding toric ideal is splittable. According to Theorem \ref{main}, this can happen if two proper subspaces of ${\rm ker}_{\mathbb{Q}}(A)$ exist
that contain all vectors corresponding to a minimal system of generators. 
\begin{Theorem} \label{Gorenstein} Let $A=\{a_{1},\ldots,a_{4}\}$ be a set of distinct positive integers, and let $C$ be the monomial curve with the
parametrization $$x_1 = t^{a_1}, x_2 = t^{a_2}, x_3 = t^{a_3}, x_4 = t^{a_4}.$$ If $C$ is a Gorenstein  monomial curve, then $I_A$ is a splittable toric ideal.
\end{Theorem}

{\em \noindent Proof.} When $C$ is a complete intersection, we derive from Theorem \ref{CompleteIntersection} that $I_A$ is splittable. Suppose that $C$ is Gorenstein and not a complete intersection monomial curve. By \cite[Theorem 3]{Bresinsky75}, $I_A$ is minimally generated by the set $$\{x_1^{d_1}- x_3^{d_{13}}x_4^{d_{14}}, x_{2}^{d_2}- x_{1}^{d_{21}}x_{3}^{d_{23}}, x_3^{d_{3}}-x_{2}^{d_{32}}x_{4}^{d_{34}}, x_{4}^{d_4}-x_{1}^{d_{41}}x_{2}^{d_{42}}, x_{1}^{d_{21}}x_{4}^{d_{34}}-x_{2}^{d_{42}}x_3^{d_{13}}\}$$ where the above binomials are unique up to isomorphism, $d_{ij}>0$, and also $$d_{1} =d_{21}+d_{41}, d_{2}= d_{32}+d_{42}, d_{3}=d_{13}+d_{23}, d_{4} =d_{14}+d_{34}.$$ Let $C=\{{\bf u}_{1}=(d_{1},0,-d_{13},-d_{14}), {\bf u}_{2}=(-d_{21},d_{2},-d_{23},0), {\bf u}_{3}=(0,-d_{32}, d_{3}, -d_{34}), {\bf u}_{4}=(-d_{41}, -d_{42}, 0, d_{4}), {\bf u}_{5}=(d_{21}, -d_{42}, -d_{13}, d_{34})\}$ and consider $C_{1}=\{{\bf u}_{1}, {\bf u}_{4}\}$ and $C_{2}=\{{\bf u}_{2}, {\bf u}_{3}, {\bf u}_{5}\}$. Since the dimension of ${\rm span}_{\mathbb{Q}}(C_{1})$ is less than or equal to 2, and the dimension of ${\rm ker}_{\mathbb{Q}}(A)$ is equal to 3, it follows that ${\rm span}_{\mathbb{Q}}(C_{1}) \subsetneqq {\rm ker}_{\mathbb{Q}}(A)$. Furthermore, since ${\bf u}_{5}=-{\bf u}_{2}-{\bf u}_{3}$, the dimension of ${\rm span}_{\mathbb{Q}}(C_{2})$ is equal to 2. However, since the dimension of ${\rm ker}_{\mathbb{Q}}(A)$ is equal to $3$, it follows that ${\rm span}_{\mathbb{Q}}(C_{2}) \subsetneqq {\rm ker}_{\mathbb{Q}}(A)$. By Theorem \ref{main}, the ideal $I_A$ is splittable. Let $A_{1}$ be the set of columns of
 \[
\left( \begin{array}{cccc}
 d_{13}d_{42} &  -d_{13}d_{41} &  d_{1}d_{42} & 0  \\
 d_{14}d_{42} & -d_{14}d_{41}-d_{1}d_{4} & 0 & d_{1}d_{42}
\end{array} \right)
	 \in \mathbb{Z}^{2 \times 4},
\] and let $A_{2}$ be the set of columns of \[
\left( \begin{array}{cccc}
 d_{2}d_{3}-d_{23}d_{32} & d_{21}d_{3} &  d_{21}d_{32} & 0  \\
 -d_{34}d_{2} & -d_{21}d_{34} & 0 & d_{21}d_{32}
\end{array} \right)
	 \in \mathbb{Z}^{2 \times 4}.
\]
Then ${\rm ker}_{\mathbb{Q}}(A_1)={\rm span}_{\mathbb{Q}}(C_1)$ and ${\rm ker}_{\mathbb{Q}}(A_2)={\rm span}_{\mathbb{Q}}(C_2)$. From the proof of Theorem \ref{main}, we have that 
$I_{A}=I_{A_1}+I_{A_2}$ is a splitting of $I_A$.
\hfill $\square$\\

In \cite[Corollary 2.13]{KT2}, it was proved that the toric ideal of a complete bipartite graph is not subgraph splittable. This means that there are no subgraphs $G_1$ and $G_2$ of $K_{m,n}$ such that $I_{K_{m,n}}=I_{G_1}+I_{G_2}$ is a splitting of $I_{K_{m,n}}$. Our aim is to prove Theorem \ref{CompleteBipartite}, which generalizes this result by proving that there are no toric ideals $I_{A_1}$ and $I_{A_2}$ such that
$I_{K_{m,n}}=I_{A_1}+I_{A_2}$ is a splitting of $I_{K_{m,n}}$.

Let $G$ be a finite, simple, connected graph on the vertex set $\{z_{1},\ldots,z_{n}\}$ with edge set $E(G)=\{e_{1},\ldots,e_{m}\}$. For each edge $e=\{z_{i},z_{j}\}$ of $G$, we associate a vector ${\bf a}_{e} \in \{0,1\}^{n}$ defined as follows: the $i$th entry of ${\bf a}_{e}$ is $1$, the $j$th entry is $1$, and all other entries are zero. We denote the toric ideal $I_{A_{G}}$ in $K[e_{1},\ldots,e_{m}]$ by $I_G$, where $A_{G}=\{{\bf a}_{e} \mid e \in E(G)\} \subset \mathbb{N}^{n}$.

Given an even closed walk $\gamma = (e_{i_1}, e_{i_2}, \ldots, e_{i_{2q}})$ of $G$, we write $B(\gamma)$ for the binomial $B(\gamma)=\prod_{k=1}^{q}e_{i_{2k-1}}-\prod_{k=1}^{q}e_{i_{2k}} \in I_{G}$. By \cite[Proposition 10.1.5]{Vil}, the ideal $I_G$ is generated by all the binomials $B(\gamma)$, where $\gamma$ is an even closed walk of $G$.

A graph $G$ is called a {\em complete bipartite} graph if its vertex set can be partitioned into two subsets $W_1$ and $W_2$ such that each edge of $G$ connects a vertex of $W_1$ to a vertex of $W_2$. It is denoted by $K_{m,n}$, where $m$ and $n$ are the numbers of vertices in $W_1$ and $W_2$, respectively. Let $W_{1}=\{v_{1},\ldots,v_{m}\}$, $W_2=\{u_{1},\ldots,u_{n}\}$ be the bipartition of the vertices of the complete bipartite graph $K_{m,n}$. By \cite[Proposition 10.6.2]{Vil}, the toric ideal $I_{K_{m,n}}$ has a minimal system of binomial generators in the form $B(\gamma)$, where $\gamma$ is a cycle of length 4. Each such cycle has four vertices $ v_a,  u_i, v_b, u_j$, where $1\leq { a}, { b}\leq m$ and $1\leq { i}, {j}\leq n$. By ${\bf 1}_{a,i}$, we denote the vector in $\mathbb{Q}^m\times \mathbb{Q}^n$ with all components zero except the one in position $a,i$, which is equal to 1. Then, by $(a,b)(i,j)$, we denote the vector ${\bf 1}_{a,i}-{\bf 1}_{b,i}+{\bf 1}_{b,j}-{\bf 1}_{a,j}$. We have two identities: $(a,b)(i,j)+(a,b)(j,l)=(a,b)(i,l)$ and $(a,b)(i,j)+(b,c)(i,j)=(a,c)(i,j)$. 

For the sake of simplicity, we shall write ${\rm ker}_{\mathbb{Q}}(K_{m,n})$ instead of ${\rm ker}_{\mathbb{Q}}(A_{K_{m,n}})$. The graph $K_{m,n}$ is bipartite with $mn$ edges and $m+n$ vertices; thus, from \cite[Proposition 10.1.20]{Vil}, the height of $I_{K_{m,n}}$ is $mn-(m+n)+1=(m-1)(n-1)$, and therefore, from \cite[Corollary 8.2.21]{Vil}, the dimension of the vector space ${\rm ker}_{\mathbb{Q}}(K_{m,n})$ is $(m-1)(n-1).$ 

The {\em complete graph} $K_{n}$ is a graph consisting of $n$ vertices, where each vertex is connected to every other vertex.

\begin{Theorem} \label{graph} Consider the complete graph $K_n$ with $n \geq 2$ vertices. Let $G_1$ and $G_2$ be two subgraphs of $K_n$ such that the union of their edge sets, $E(G_1) \cup E(G_2)$, equals the edge set of $K_n$, $E(K_n)$. Then, at least one of these subgraphs contains a spanning tree of $K_n$.
\end{Theorem}
{\em \noindent Proof.} The proof of the theorem will be conducted by induction on $n$. The case $n=2$ is straightforward and evident. Assume the proposition holds for $n=l$. Let $G_1$ and $G_2$ be two subgraphs of $K_{l+1}$ such that $E(G_1)\cup E(G_2)=E(K_{l+1}).$
Let $G$ be an induced subgraph of $K_{l+1}$ with $l$ vertices. Then the graph $G$ is complete. We have $E(G\cap G_1)\cup E(G\cap G_2)=E(G)$. According to the hypothesis, one of the graphs $G \cap G_1$ or $G \cap G_2$ contains a spanning tree of $G$. Assume it is $G \cap G_1$. Let $v$ be the vertex of $K_{l+1}$ that does not belong to $G$. We distinguish the following cases.
\begin{enumerate}
\item There exists an edge of $G_1$ that has $v$ as one of its vertices. This edge, together with the spanning tree of $G$, forms a spanning tree of $K_{l+1}$ that is contained within $G_1$.
\item None of the edges in $G_1$ have $v$ as one of their vertices. Since $E(G_1)\cup E(G_2)=E(K_{l+1})$, each edge of $K_{l+1}$ containing $v$ is also an edge of $G_2$. Consequently, all edges of $K_{l+1}$ containing $v$ form a spanning tree of $K_{l+1}$ that is contained within $G_2$. \hfill $\square$\\
\end{enumerate}

\begin{Theorem} \label{CompleteBipartite} The toric ideal of a complete bipartite graph $K_{m,n}$ is not splittable. 
\end{Theorem}
{\em \noindent Proof.} 
Suppose that the toric ideal $I_{K_{m,n}}$ is splittable. According to \cite[Theorem 2.3]{HO}, the toric ideal $I_{K_{m,n}}$ has a unique minimal system of binomial generators in the form 
$\{B({\bf u}) \mid {\bf u}\in C\}$, where $$C=\{(x,y)(s,t) \mid (v_{x},u_{s},v_{y},u_{t}) \ \textrm{is a cycle of} \ K_{m,n}\}.$$
It follows that ${\rm span}_{\mathbb{Q}}(C) = {\rm ker}_{\mathbb{Q}}(K_{m,n})$.

By Theorem \ref{main}, there exist 
sets $C_{1}$ and $C_2$ such that $C=C_1 \cup C_2$, ${\rm span}_{\mathbb{Q}}(C_1)\subsetneqq {\rm ker}_{\mathbb{Q}}(K_{m,n})$, and ${\rm span}_{\mathbb{Q}}(C_2)\subsetneqq {\rm ker}_{\mathbb{Q}}(K_{m,n})$. Fix two vertices, $v_{x}$ and $v_{y}$, of $W_1$, and consider all vectors of the form $(x,y)(s,t)$, where $(v_{x},u_{s},v_{y},u_{t}) \ \textrm{is a cycle of} \ K_{m,n}$. Let $V^{(x,y)}={\rm span}_{\mathbb{Q}}((x,y)(s,t) \mid u_{s}, u_{t} \in W_{2})$ and $G^{(x,y)}_i$ subgraphs of $K_n$ such that  $E(G^{(x,y)}_i)=\{\{u_{s},u_{t}\} \mid (x,y)(s,t)\in C_i\}\subset K_n$, for $i=1$ or $2$. Then, $E(G^{(x,y)}_1) \cup E(G^{(x,y)}_2) = E(K_n)$, and by Theorem \ref{graph}, we have that at least one of the subgraphs $G^{(x,y)}_1$ or $G^{(x,y)}_2$ contains a spanning tree of $K_n$. Suppose that $G^{(x,y)}_i$ contains a spanning tree $T$ of $K_n$. Since $T$ is a spanning tree, there exists a path $(b_1, \ldots, b_r)$ in $T$ such that $u_{s}=u_{b_1}$ and $u_{t}=u_{b_r}$. Then, $(x,y)(s,t)=(x,y)(b_1,b_r)=(x,y)(b_1,b_2)+(x,y)(b_2,b_3)+\cdots+(x,y)(b_{r-1},b_r)$. Note that by the construction of the graph $G^{(x,y)}_i$, all vectors $(x,y)(b_j,b_{j+1})$, $1 \leq j \leq r-1$, belong to $C_i$. Hence, we have $V^{(x,y)}\subset {\rm span}_{\mathbb{Q}}(C_i)$. Therefore, for each $v_{x}$ and $v_{y}$ in $W_1$, either $V^{(x,y)}\subset {\rm span}_{\mathbb{Q}}(C_1)$ or $V^{(x,y)}\subset {\rm span}_{\mathbb{Q}}(C_2)$. We define two subgraphs of $K_m$ by $E(G_i)=\{\{v_{x},v_{y}\}|V^{(x,y)}\subset {\rm span}_{\mathbb{Q}}(C_i)\}$ for $i=1$ or $2$. Then, we have $E(G_1) \cup E(G_2) = E(K_m)$, and by Theorem \ref{graph}, we have that at least one of these, say $G_1$, contains a spanning tree $S$ of $K_m$. Consider the vector $(x,y)(s,t)$. Since $S$ is a spanning tree  there exists a path $(a_1, \ldots, a_k)$ in $S$ such that $v_{x}=v_{a_1}$
and $v_{y}=v_{a_k}$. Then, $(x,y)(s,t)=(a_1,a_k)(s,t)=(a_1,a_2)(s,t)+(a_2,a_3)(s,t)+\cdots+(a_{k-1},a_k)(s,t)$. Note that, by the construction of the graph $G_1$, all vectors $(a_j,a_{j+1})(s,t)$ , $1 \leq j \leq k-1$, belong to $ {\rm span}_{\mathbb{Q}}(C_1)$. Therefore, it follows that ${\rm ker}_{\mathbb{Q}}(K_{m,n}) = {\rm span}_{\mathbb{Q}}(C) \subset {\rm span}_{\mathbb{Q}}(C_1)\subset {\rm ker}_{\mathbb{Q}}(K_{m,n})$. This implies that ${\rm span}_{\mathbb{Q}}(C_1)= {\rm ker}_{\mathbb{Q}}(K_{m,n})$, which is a contradiction. Therefore, the toric ideal $I_{K_{m,n}}$ is not splittable.
\hfill $\square$

\end{document}